\documentstyle{article}

\title{Generic Integral Manifolds for Weight Two Period Domains}
\author{James A. Carlson and Domingo Toledo}
\date{June 2, 2001}

\newtheorem{theorem}{Theorem}
\newtheorem{lemma}{Lemma}
\newtheorem{proposition}{Proposition}
\newtheorem{definition}{Definition}

\newfam\meuffam
\font\tenmeuf=eufm10
\font\sevenmeuf=eufm7
\font\fivemeuf=eufm5

\textfont\meuffam=\tenmeuf
\scriptfont\meuffam=\sevenmeuf
\scriptscriptfont\meuffam=\fivemeuf
\def\germ{\fam\meuffam\tenmeuf}

\newcommand{\GG}{{\cal G}}
\newcommand{\gu}{{\germ u}}
\newcommand{\ga}{{\germ a}}
\newcommand{\half}{{1 \over 2 }}
\newcommand\set[1]{\{\ #1 \ \}}
\newcommand\sett[2]{\{\ #1 \ \ |\ #2 \ \}}

\newcommand\pd[2]{{ \partial {#1} \over \partial {#2} }}
\newcommand{\mmpd}[3]{{\partial^2 #1 \over \partial #2 \partial #3 }}
\newcommand{\mmmpd}[4]{{\partial^3 #1 \over \partial #2 \partial #3 \partial #4 }}

\newcommand{\commadots}{\,,\,\cdots\,,\,}
\newcommand{\map}{\longrightarrow}

\newcommand{\R}{{\bf R}}
\newcommand{\C}{{\bf C}}
\newcommand{\eqdef}{\ \ \smash{ 
   \mathop{=}\limits_{def}}\ \ } 
\newcommand{\tr}{{}^t\kern-0.9pt} 
\newcommand{\Tr}{{}^T\kern-0.9pt} 

\begin{document}

\maketitle

\begin{abstract}
We define the notion of a generic integral element for the
Griffiths distribution on a weight two period domain, 
draw the analogy with the classical contact distribution, and
then show how to explicitly construct an infinite-dimensional
family of integral manifolds tangent to a given element.

\end{abstract}

\section{Introduction}

In this note we shall study a class of integral manifolds for a
generalization of the classical contact distribution.
This distribution arises naturally in algebraic geometry as
the infinitesimal restriction satisfied by the period maps
of algebraic surfaces.  To be more precise, consider the
group $G = SO(2p,q)$, the compact subgroup $U(p)\times SO(q)$,
and the homogeneous space
\[
    D = G / V.
\]
The homogeneous space parametrizes Hodge decompositions
\[
  H_{\C} = H^{2,0} \oplus H^{1,1} \oplus H^{0,2}
\]
on a fixed complex vector space endowed
with a nondegenerate, symmetric bilinear form
defined over an underlying real vector space
$H_\R$ whose complexification is $H_\C$.  The models
for these spaces are the primitive second cohomology
of an algebraic surface with real or complex coordinates.  
If one has a family of algebraic surfaces $S_t$ parametrized
by $t$, then one has a family of Hodge decompositions
which are determined by the subspaces $H^{2,0}(S_t)$.
As Griffiths showed in \cite{Griff}, these spaces satisfy
\[
  { d \over dt } H^{2,0}(S_t) 
    \subset H^{2,0}(S_t) \oplus H^{1,1}(S_t). 
\]
This condition, which asserts that a famly of subspaces
$H^{2,0}_t$ cannot vary arbitrarily, but rather like a Frenet 
frame of a space curve, defines the
\emph{Griffiths distribution} on the homogeneous space $D$.
Thus families of algebraic surfaces define integral
manifolds of this distribution.

In \cite{CTAMS} it was shown that for $p > 1$ even and $q > 1$, 
integral manifolds which are of maximal dimension are 
quite special: as germs they are congruent under
the action of $G$; moreover, they correspond to the germ
of an imbedding $U(p,q/2) \map SO(2p,q)$.  Thus
germs of integral manifolds of this kind are rigid and depend on
finitely many parameters.

Here we shall study certain integral manifolds which are maximal 
with respect to inclusion: they do not lie in a larger dimensional
integral manifold.  There is natural class of these which we call
``generic,'' characterized by the existence of
a family of vectors $v(t)$ in $H^{2,0}_t$ 
such that the  partial derivatives
$\partial v /\partial t_i$
span $H^{1,1}_t$.
(See definition \ref{def:generic}). These
integral manifolds behave quite differently: they are not rigid, and 
they depend on infinitely
many parameters.  Moreover, as we shall see in Theorem \ref{maintheorem}
and the  paragraph which follows it, generic integral
manifolds are given quite explicitly by a system of generating functions.
This is in strict analogy to the case of the contact distribution,
where maximal integral manifolds are flexible and determined
explicity by a generating function.  The function(s) constitute
the infinite-dimensional parameters.

\section{Generalized Contact Distributions}

Recall that the contact distribution is the annihilator $E$ in the
tangent bundle of
$\R^{2n+1}$ or $\C^{2n+1}$ of the one-form
$$
    \omega = dz - y\cdot dx = dz - \sum_{i=1}^n y_idx_i ,
    \label{contactform}
$$
and that integral manifolds are by definition submanifolds
tangent to $E$. They are of dimension at
most $n$ and those of maximum dimension --- the Legendre manifolds
(see \cite{Arnold}) --- are given, up to certain admissible
changes of coordinates, by a generating function
$f$ through the prescription
\begin{eqnarray}
  & z = f(x) \\
  & y = \nabla f .
\end{eqnarray}
Thus maximal integral manifolds are graphs
of one-jets of functions, and so constitute an infinite-dimensional
family.

After a
suitable change of variables, any distribution of codimension 
one can be reduced locally  to a product of the trivial
distribution and a contact distribution.  Consequently the local
nature of their maximal integral manifolds
is completely understood.  By ``maximal'' we mean ``maximal with
respect to inclusion.''  The situation for distributions of
codimension greater than one is completely different.  There is no
general theory which answers the basic questions and the nature of
the integral manifolds is in general very complicated.  There is
always a maximum dimension for integral manifolds and manifolds
of that dimension are obviously maximal with respect to inclusion. 
However, maximal integral manifolds may not be of maximal dimension.
The distribution which we shall study here is of codimension greater
than one and exhibits the just-mentioned behavior.  It is a natural
one to study for a number of reasons. First, it is defined by a
matrix-valued analogue of the contact distribution.  Second, it
arises as a local model $U$ for the Griffiths distribution
on a period domain $D$ of weight two \cite{Griff} discussed in the
introduction. 

We shall now describe this local model $U$, then give the main result
of the paper. To this end, consider the group
$G$ of matrices of the form
$$
  g = \left(\matrix{
                  1_p & 0 & 0 \cr
                  X & 1_q & 0 \cr
                  Z & Y & 1_p \cr
                    }\right)
$$
where $1_n$ denotes the $n\times n$ identity matrix.
A basis for the left-invariant one-forms on $G$ is given by the
Maurer-Cartan form
\begin{equation}
  \Omega \eqdef 
         g^{-1}dg = \left(\matrix{
                  0 & 0 & 0 \cr
                  dX & 0 & 0 \cr
                  \omega & dY & 0 \cr
                    }\right)
 \label{MCform}
\end{equation}
where
\begin{equation}
  \omega = dZ - YdX
 \label{matrixcontactform}
\end{equation}
On $G$ define a distribution $E$ as the set of tangent vectors
which annihilate the entries of $\omega$.  When $p = 1$, the group is
the Heisenberg group and $\omega$ is the contact form.  Our local
model $U$ is the unipotent subgroup defined by the equations
\begin{equation}
    \quad Y  = \tr X
\label{pereq-a}
\end{equation}
and
\begin{equation}
    \quad Z + \tr Z = \tr X X .
 \label{pereq-b}
\end{equation}
A neighborhood of the identity in $U$ is
isomorphic to a neighborhood of a fixed but arbitrary point of $D$.
Under this identification the restriction of the Griffiths distribution
is the same as the restriction of the distribution $E$ to $U$. 

Since $X$ determines $Y$ and the
symmetric part of $Z$, coordinates on $U$ are given by the entries of
$X$ and the skew-symmetric part of $Z$. Therefore $U$ has dimension
$pq + p(p-1)/2$.   From (\ref{matrixcontactform}) and
(\ref{pereq-a}), one finds that
\begin{equation}
    \omega = dZ  - \tr X dX,
    \label{Ucontacteq-a}
\end{equation}
and from the exterior derivative of (\ref{pereq-b}) one finds that
\begin{equation}
     \omega^+ = 0,
\label{Ucontacteq-b}
\end{equation}
where $\omega = \omega^+ + \omega^-$ is the decomposition into
symmetric and antisymmetric parts.  Thus $E$ has dimension $pq$ and
codimension $p(p-1)/2$.  Our main interest will be in the case $p >
2$, i.e., codimension greater than one.

Consider now a subspace $S$ of the tangent space to $U$ at some point.
If it is tangent to an integral manifold it annihilates not only
$\omega$ but also $d\omega$.  An arbitrary subspace satisfying these
two conditions --- a potential tangent space to an integral manifold
--- is called an {\sl integral element}.  For the contact distribution
all integral elements are integrable.  For the Griffiths distribution
``integrability'' holds for integral elements of maximal dimension
($pq/2$ when $q$ even, $p(q-1)/2 + 1$ when $q$ odd, \cite{CKT}).  For other integral
elements, e.g., those which are maximal with respect to inclusion
but not of maximal dimension, one could presumably answer the integrability
question using the  Cartan-K\"ahler theory.  What we do instead
is to solve the integrability problem explictly for generic integral
elements, which are easily shown to be maximal \cite{CTAMS}:

\begin{theorem} Let $S$ be a generic $q$-dimensional integral
element for a period domain with Hodge numbers $p = h^{2,0}$, $q = h^{1,1}$,
where $p > 1$. Then $S$ is tangent to an integral
manifold. Such integral manifolds are determined in a canonical way by holomorphic functions
$f_2 \commadots f_p$ of a complex variable $u = (u_1 \commadots
u_q)$ which satisfy the system of partial differential equations
\begin{equation}
   [H_{f_i}, H_{f_j}] = 0
    \label{Hcommutatoreq}
\end{equation}
where $H_f$ is the Hessian matrix of $f$.  The space of solutions to
this equation for fixed $S$ is infinite-dimensional.
\label{maintheorem}
\end{theorem}

\noindent
The set of generic integral elements is open in the set of all
$q$-dimensional integral elements.  We shall formulate
and prove this result in the next section.  To explain the canonical
construction, recall that the entries of the matrices $X$ and the
skew-symmetric part of $Z$ give coordinates on $U$.  Thus an integral
manifold will be specified by giving
these coordinates in terms of the functions $f_i$. To do
so, let $[a_1 \commadots a_p]$ denote the matrix with
column vectors
$a_i$, set
$$
  X(u) = [\, u, \nabla f_2 \commadots \nabla f_p \, ],
$$
and put
$$
  Z_{j1}(u) = f_j(u),
$$
where $j > 1$.  For the entries $Z_{jk}$ with $j>k$, choose arbitrary
solutions of the equations
$$
   dZ_{jk} = \sum_a X_{aj}dX_{ak}
$$
deduced from the $jk$ entry of 
$$
  \omega = dZ - \tr X dX = 0.
$$  
Use these to
determine the antisymmetric part of $Z$.  For the symmetric part use
the quadratic equation (\ref{pereq-b}).  The analogy with the
contact system --- both in the form of the equations and the form of
the solutions, is clear.  Indeed, when $p = 2$ the equation $\omega =
0$ is the same as the equation
$$
   dZ_{12} - \sum_a X_{a1}dX_{a2} = 0
$$
Thus, if we set $z = Z_{12}$, $x = (x_{11} \commadots x_{1q})$.
and $y = (x_{21} \commadots x_{2q})$, then both equations and
solutions coincide with those of the contact case.  Note, however,
that the relations (\ref{Hcommutatoreq}) are a new feature of the
case $p > 2$.

\section{Integral elements}

In order to give a precise definition of ``generic'' we describe in
some detail the tangent space of $D$ at a fixed point of reference
and the integral elements it contains.  To this end we choose
the local correspondence between
$D$ and $U$ so that identity matrix of $U$ is mapped to the reference
point.  Thus the Lie algebra $\gu$ --- the tangent space at the
identity of $U$ --- corresponds to the tangent space of $D$ at the
reference point.  Now consider a curve $g(t)$ based at the identity
matrix with arbitrary initial tangent $\tau = g'(0)$.
It has the form
$$
  g(t) = \left(\matrix{
                  1_p & 0 & 0 \cr
                  a(t) & 1_q & 0 \cr
                  b(t) & \tr a(t) & 1_p \cr
                    }\right)
$$
where 
$$
b = c(t) + \half \tr a(t)a(t)
$$
with $c(t)$ is skew-symmetric, and where $a(0) = 0$, $c(0)
= 0$.  Differentiating, we find
$$
  \tau = \left(\matrix{
                  0 & 0 & 0 \cr
                  \phi & 0 & 0 \cr
                  \psi & \tr \phi & 0 \cr
                    }\right)
$$
where $a'(0) = \phi$ and $c'(0) = \psi$ are arbitrary matrices subject
to the condition that $\psi$ be skew-symmetric.  We can read this as
saying $\Omega(\tau) = \tau$, where $\Omega$ is the Maurer-Cartan form
(\ref{MCform}).  In more detail,
\begin{eqnarray}
  & dX(\tau) = \phi \\
  & dZ(\tau) = \psi .
\end{eqnarray}
Thus a matrix $\tau(\phi, \psi)$ is an element of $E$ if an only if
$\psi = 0$.  Consequently the map
$$
   \phi \mapsto \tau(\phi) \eqdef\tau(\phi,0)
$$
defines an isomorphism of $E$  at the identity with $q \times p$
matrices, i.e., with linear maps
$$
  \phi:\C^p \map \C^q
$$

Now let $\ga \subset \gu$ be an integral element and let
$\tau_i = \tau(\phi_i)$ be vectors in $\ga$.  By definition
$\ga$ annihilates
$$
   d\omega = -\tr dX \wedge dX.
$$
But
\begin{eqnarray}
   d\omega(\tau_1, \tau_2) 
    &= 
   - \tr  dX(\tau_1) dX(\tau_2) + \tr  dX(\tau_2) dX(\tau_1) \cr
    &=
   - \tr \phi_1 \phi_2 + \tr \phi_2 \phi_1 \cr
\end{eqnarray}
so that the commutator
\begin{equation}
  (\phi_1, \phi_2) 
    \eqdef
  \tr \phi_1 \phi_2 - \tr \phi_2 \phi_1
  \label{newcommutator}
\end{equation}
vanishes.  Equivalently, the Lie bracket $[\tau_1, \tau_2]$ vanishes.
Thus we may regard $\ga$, either as a subspace of $\gu$ or of the linear
maps from $\C^p$ to $\C^q$, as an {\sl abelian subspace}.  To
summarize:

\begin{lemma} A subspace of $\gu$ is an integral element if and only
if it is an abelian subspace.
\end{lemma}
  
\noindent
We can now define what me mean by generic:

\begin{definition} A $q$-dimensional abelian subspace $\ga$ of $\gu$
is {\sl generic} if there is a vector $v \in \C^p$ such that
$\ga(v) = \C^q$.
\label{def:generic}
\end{definition} 

\noindent
By $\ga(v)$ we mean the space $\sett{ \phi(v) }{\phi\in\ga }$.  To
justify the terminology we claim that (a) there are such spaces and
(b) the condition that a space be generic is an open one.  For the
first point let $\set{ e_i }$ denote the standard basis for $\C^n$
and let $v \cdot w$ denote the complex dot product.  Then the
commutator as defined in (\ref{newcommutator}) of $q \times p$
matrices
$a = [a_1 \commadots a_p]$ and $b = [b_1 \commadots b_p]$,  
is the skew-symmetric matrix with entries
$$
  (a,b)_{ij} = a_i \cdot b_j - b_i \cdot a_j .
$$
Set
$$
   M_i = [e_i, 0 \commadots 0 ]
$$
for $i = 1..q$ and let $\ga_0$ be their span.  It is clear that
$(M_i,M_j) = 0$, so that $\ga_0$ is abelian.  For the second point
consider the space $A$ of framed abelian subspaces, that is, abelian
subspaces endowed with a basis $\set{ M_i }$. Consider the function
$F_v$ on $A$ defined by
$$
   (M_1 \commadots M_q) 
      \mapsto 
    M_1(v)\wedge \cdots \wedge M_q(v)
$$
The union $\GG$ of the sets $\GG(v) = \sett{ (M_i) = \hbox{basis for
an abelian space}}{ F_v(M_i)\ne 0 }$ is open in $A$ and contains any
framing of $\ga_0$, whence the claim. 

Now consider the $f$ transformation defined by
\begin{eqnarray}
  & X \map BXA \cr
  & Z \map \tr AZA \cr
\end{eqnarray}
where $X$ is a $q\times p$ matrix, $Z$ is a $p\times p$ matrix,
where $A$ is invertible and where $B$ is complex orthogonal.  
The set of such transformations constitutes a complex Lie group
$H \cong GL(p,\C)\times O(p,\C)$ which acts on the local model $U$. 
This action fixes the identity and acts on the form $\omega$ by
$$
  \omega \map \tr A \omega A
$$
Consequently it preserves the distribution $E$ and so maps integral
manfifolds to integral manifolds.  Therefore in studying integral
manifolds of $E$ we may do so up to the action of $H$. Note also that
for commutators,
$$
  (BX_1A, BX_2A) = \tr A(X_1,X_2)A,
$$
so that the transform of an abelian space is an abelian space. 
From this one sees that the orbit of $\GG(e_1)$ is
$\GG$;  consequently, we may, without loss of generality, reason about
$\GG(e_1)$ in place of $\GG$.  But an element of $\GG(e_1)$ is a
$q$-dimensional abelian space with a basis elements of the form
$$
  M_j = [\, e_j, * \commadots * \,],
$$
where ``$*$'' stands for a column vector.  We shall call
such bases ``distinguished.''

Distinguished bases for abelian spaces can be characterized as
follows.  Given a matrix $A$, let $(A)_k$ denote the $k$-th column. 
Consider next a system of $q\times q$ matrices  $\set{A_j}$, 
where $j = 2 \dots p$, and  construct a new system
of $q\times p$ matrices
$$
   M_k = [ e_k,\ (A_2)_k \commadots (A_p)_k  ],
$$
where $k = 1 \dots q$.
The correspondence $\set{ A_j } \map \set{ M_k} $
is one-to-one, and a routine computation shows the
following:

\begin{proposition} The span of a distinguished basis $\set{M_k}$, 
$k = 1 \dots q$,
is an abelian space if and only if $\set{ A_j }$, $j = 2 \dots p$,
is a commuting set of
symmetric matrices.  The span of the $M_k$ is then a maximal abelian space.
\label{distbasisprop}
\end{proposition}

\noindent
Moreover, it is not hard to show that the natural relation between
Hessians and tangent spaces holds:

\begin{proposition} Let $\set{ A_j }$ be a commuting set of symmetric
matrices and let $\set{ f_j }$ be a solution to (\ref{Hcommutatoreq})
such that $f_j(0) = 0$, $\nabla f_i(0) = 0$, and $H_{f_j}(0)  = A_j$. 
Then the tangent space at the identity to the associated integral
manifold is the abelian space associated to $\set{ A_j }$ which 
has basis $\set{ M_k }$.
\end{proposition}

\section{The canonical construction}

We will now show that an integral manifold of $E$
whose tangent space at the
identity is in $\GG(e_1)$ is given locally by the canonical construction. To this
end, note that the genericity hypothesis (the condition $\ga(e_1) = \C^q$) is
equivalent to the condition that the components of $(dX)_1$ are independent.
In this case 
$$
   \eta = dX_{11} \wedge \cdots \wedge dX_{r1}
$$
is nonzero.  By shrinking $M$ we may assume that the product 
$\eta$ is nonzero on all of
$M$, so that the functions $X_{a1}$ are independent on it.  Then
(\ref{Ucontacteq-a}) implies that
$$
  dZ_{ij} = \sum X_{ai}dX_{aj} .
$$
Consider in particular the case $j = 1$, $i > 1$, for which we obtain
the equation
$$
  dZ_{i1} = \sum X_{ai}dX_{a1}.
$$
Since no $dX_{ab}$ for $b \ne 1$ occur, we see that $Z_{i1}$ may be
viewed as function $f_i$ of the entries of $(X)_1$, the first column of
$X$. Moreover, the $X_{ai}$ are functions of these same entries,
namely,
$$
    X_{ai} = \pd{f_i}{X_{a1}}
$$
So far we have used just some of the equations $\omega = 0$ determined by
(\ref{Ucontacteq-a}).  These equations assert that there exist certain
functions 
$Z_{ij}$ for $i > j > 1$ of  $(X)_j$ which are in turn
functions of $(X)_1$.  Consequently the already-determined forms
$\phi_{ij} = dZ_{ij}$ must be closed.  Now
$$
  \phi_{ij} 
    = 
     \sum_a 
       \pd{f_i}{X_{a1}} d\left(\pd{f_j}{X_{a1}}\right)
    =
      \sum_{ab} 
        \pd{f_i}{X_{a1}} \mmpd{f_j}{X_{a1}}{X_{b1}}{ dX_{b1}}
$$
and so
$$
  d\phi_{ij}
    = 
       \sum_{abc} \left[ 
         \mmpd{f_i}{X_{a1}}{X_{c1}} \mmpd{f_j}{X_{a1}}{X_{b1}} dX_{c1}\wedge dX_{b1}
           +
          \pd{f_i}{X_{a1}} \mmmpd{f_j}{X_{a1}}{X_{b1}}{X_{c1}} dX_{c1}\wedge dX_{b1}
          \right]
$$
The third partial derivative is symmetric in $b$ and $c$ whereas the product
$dX_{c1}\wedge dX_{b1}$ is antisymmetic in these indices.  Consequently the 
second sum vanishes.  Therefore the first sum must vanish.  Since the coefficient
of $dX_{c1}\wedge dX_{b1}$ is symmetric in $b$ and $c$ we conclude that
it must vanish.  But inspection reveals that coeffient to be the $bc$ entry
of the commutator
$$
  [H_{f_i},H_{f_j}],
$$
where $H$ denotes the Hessian matrix.  Consequently the consistency of our overdetermined
system is just the set of partial differential equations
(\ref{Hcommutatoreq}).

\section{Existence results}

Let us now consider the problem of the existence and nature of solutions
to the system of partial differential equations $[H_{f_i},H_{f_j}] = 0$.
The most basic question is whether enough there are enough solutions to
pass an integral manifold through an arbitrary integral element.  By Proposition
(\ref{distbasisprop}) this is
equivalent to the problem of constructing functions whose Hessians commute
and whose values at the origin are given commuting symmetric matrices $\set{A_\ell}$.
For these it is enough to take the quadratic functions
$$
   f_\ell = \half\sum_{ij} (A_\ell)_{ij} u_iu_j
$$
We have therefore shown the following and with it part of (\ref{maintheorem}):

\begin{proposition} Any generic element is tangent to an integral manifold.
\label{integrabilityprop}
\end{proposition}

Let us consider now the problem of finding additional solutions to 
the equations (\ref{Hcommutatoreq}). Note first that
if two functions $f$ and $g$ have diagonal Hessians then they
automatically satisfy $[H_f,H_g] = 0$.
The condition that the Hessians be diagonal is the 
overdetermined system of equations
$$
   \mmpd g{x_i}{x_j} = 0 \qquad \hbox{for $i \ne j$},
   \eqno{(*)}
$$
In the case of two variables there is just one equation, a form of the wave equation,
which has solutions of the form $h_1(x_1) + h_2(x_2)$.  The solutions in the 
general case have the same form,
$$
  g(x_1 \commadots x_n) = \sum_i h_i(x_i),
  \label{diagsolution}
$$
where the $h_i$ are arbitrary functions of one variable.
Therefore an integral manifold is specified by a set of functions
$\set{ h_{ij}(u) }$, where
$$
  f_i(x_1 \commadots x_s) = \sum_j h_{ij}(x_j) .
$$ 
This does not give a complete set of solutions, but it does give an infinite-dimensional
set.  Moreover, we can choose the $h_{ij}$ in such a way that so that the Hessian of 
$f_i$ at the origin is an arbitrary diagonal quadratic form.  In fact, we can do somewhat
more.  Let $\set{ A_\ell }$ be a set of commuting complex symmetric matrices at least one
of which has distinct eigenvalues.  Then there is a set of common eigenvectors which form
a basis for $\C^q$ and which are orthogonal relative to the complex dot product. 
Consequently there is a complex orthogonal matrix $C$ which simultaneously diagonalizes
the $A_\ell$.  Let $D_\ell$ be the diagonal matrix corresponding to $A_\ell$, where 
$CA_\ell \tr C = D_\ell$.  There is an infinite dimensional family of solutions $(f'_\ell)$
to (\ref{Hcommutatoreq}) such that $H_{f'_\ell}(0) = D_\ell$.  Let $f_\ell(x) = f'_\ell(Cx)$.
Then
$$
   H_{f_\ell} = \tr C H_{f'_\ell} C
$$
and, since $C$ is complex orthogonal,
$$
   [H_{f_\ell},H_{f_m}]  = \tr C [ H_{f'_\ell}, H_{f'_m}] C
$$
Thefore the infinite-dimensional family of functions $(f_\ell)$ is also a set of 
solutions to (\ref{Hcommutatoreq}), and each member has the specified initial Hessians
$(A_\ell)$. Consequently our previous integration result (\ref{diagsolution}) can be
strengthened:

\begin{theorem} Any generic abelian space of dimension $r$ is tangent to an integral
manifold.  Moroever,  the set of germs of integral manifolds tangent to this space is
infinite dimensional.
\end{theorem}

\subsection*{Remark.} Consider the case $q = 2$ with $p > 2$ arbitrary.
Fix $f_2$ arbitrarily but generically in the sense that the Hessian
generically has distinct eigenvalues.  Consider the equations 
$[H_{f_2},H_{f_j}] = 0$, for $j > 2$.  For given $j$ one has a single
non-trivial partial differential equation which is linear of second
order in $f_2$.  By the Cauchy-Kowaleska Theorem \cite{FJohn} the solution
space is infinite-dimensional. Now let $V(u)$ be the matrix of eigenvectors
of $H_{f_2}(u)$.  Then the transformation $A \map \tr V(u)AV(u)$ simultaneously
diagonalizes all of the matrices $H_{f_i}(u)$.  Consequently all of these
matrices commute with each other, i.e., the functions $f_i$ solve (\ref{Hcommutatoreq}).  
The case $q > 2$ is more
complicated because the  partial differential equations constitute an overdetermined system, somewhat
like the system $\nabla f  = \xi$ for a given vector field $\xi$.

\end{document}